\documentclass{article}

\usepackage{latexsym}
\usepackage{amssymb}
\usepackage{amsfonts}
\usepackage{amsmath}
\usepackage[all]{xy}

\newtheorem{thm}{\bf Theorem}[section]
\newtheorem{cor}[thm]{\bf Corollary}

\newtheorem{lem}[thm]{\bf Lemma}
\newtheorem{prop}[thm]{\bf Proposition}
\newtheorem{defn}[thm]{\bf Definition}

\newtheorem{rem}[thm]{\bf Remark}

\newtheorem{exmps}[thm]{\bf Examples}

\def\proof{{\parindent0pt {\bf Proof.\ }}}

\def\fd{{\rm fd}}

\def\GFD{{\rm  GFD}}

\def\lGCFD{{  {\rm l.G}_{\rm C}-{\rm FD}}}

\def\GCFD{{\rm  G_C\,\textnormal{-}\,FD}}

\def\FPid{{\rm  FP-id}}

\def\fc{{\rm \mathcal{F}_C\,\textnormal{-}\,fd}}
\def\ic{{\rm \mathcal{I}_{C^+}\,\textnormal{-}\,id}}

\def\F{\mathcal{F}_C}

\def\I{\mathcal{I}_{C^+}}
\def\C{\mathcal{C}_C}

\def\Y{\mathcal{Y}}
\def\X{\mathcal{X}}
\def\G{\mathcal{G}}
\def\H{\mathcal{H}}
\def\V{\mathcal{V}}

\def\XX{\textbf{X}}

\def\A{\mathcal{A}}
\def\B{\mathcal{B}}
\def\Q{\mathcal{Q}}
\def\R{\mathcal{R}}
\def\W{\mathcal{W}}
\def\M{\mathcal{M}}

\def\Gfd{{\rm Gfd}}

\def\Flat{\mathcal{F}}

\def\H{\mathcal{H}}
\def\V{\mathcal{V}}

\def\Gfc{{  {\rm G}_{\rm C}\,\textnormal{-}\,{\rm fd}}}

\def\GC{{\rm G}_C}
\def\GCF{{\rm G}_{\rm C}{\rm F}}
\def\GCC{{\rm G}_{\rm C}{\rm C}}

\def\Im{{\rm Im}}

\def\Ker{{\rm Ker}}
\def\Ho{{\rm Ho}}

\def\resdim{{\rm resdim}}
\def\coresdim{{\rm coresdim}}

\def\Ext{{\rm Ext}}
\def\Tor{{\rm Tor}}
\def\Hom{{\rm Hom}}
\def\End{{\rm End}}
\def\Prod{{\rm Prod}}
\def\Add{{\rm Add}}

\def\Inj{\mathcal{I}}
\def\Flat{\mathcal{F}}
\def\Cot{\mathcal{C}}
\def\Proj{\mathcal{P}}

\def\sup{{\rm sup}}

\newcommand{\GF}{\rm GF}

\newcommand{\cqfd}
{\rule{2mm}{2mm}%
\medbreak%
\par%
}

\begin{document}
\title{Relative weak global Gorenstein dimension, AB-contexts and model structures}

\author{Driss Bennis$^1$ \hskip 2cm  Rachid El Maaouy$^{2,}$\footnote{Corresponding author} \\ \\ J. R. Garc\'{\i}a Rozas$^3$ \hskip 1,5cm Luis Oyonarte$^4$}

\date{}

\maketitle

{\small1: CeReMaR Research Center, Faculty of Sciences, B.P. 1014, Mohammed V University in Rabat, Rabat, Morocco.
	
	\noindent e-mail address: driss.bennis@um5.ac.ma; driss$\_$bennis@hotmail.com
	
	2: CeReMaR Research Center, Faculty of Sciences, B.P. 1014, Mohammed V University in Rabat, Rabat, Morocco.
	
	\noindent e-mail address: rachid\_elmaaouy@um5.ac.ma; elmaaouy.rachid@gmail.com
	
	3: Departamento de  Matem\'{a}ticas, Universidad de Almer\'{i}a, 04071 Almer\'{i}a, Spain.
	
	\noindent e-mail address: jrgrozas@ual.es
	
	4: Departamento de  Matem\'{a}ticas, Universidad de Almer\'{i}a, 04071 Almer\'{i}a, Spain.
	
	\noindent e-mail address: oyonarte@ual.es}
	
	\begin{abstract}
	
	In this paper we introduce and study the weak Gorenstein global dimension of a ring $R$ with respect to a left $R$-module $C$. We provide several characterizations of when this homological invariant is bounded. Two main applications are given: first, we prove that the weak Gorenstein global dimension of $R$ relative to a semidualizing $(R,S)$-bimodule $C$  can be computed either by the $\GC$-flat dimension of the left $R$-modules or right $S$-modules, just like the (absolute) weak global dimension. As a consequence, a new argument for solving Bennis' conjecture is obtained.  As a second application, we give a concrete description of the weak equivalences in the $\GC$-flat model structure recently found by the authors. In order to prove this result, an interesting connection between abelian model structures and (weak) AB-weak contexts is proved. This connection leads to a result that can be applied to obtain abelian model structures with a simpler description of trivial objects.
	
\end{abstract}


	
	\normalsize
\section{Introduction}

 \hskip .5cm Auslander, in 1966, introduced the notion of G-dimension of a finitely generated module over a Cohen-Macaulay noetherian ring and studied their basic properties. This notion of G-dimension has been extended in several directions. For instance, Enochs et al. (\cite{EJT94}, \cite{EJ95}) introduced the Gorenstein projective dimensions and the Gorenstein flat dimensions. Golod \cite{Gol84} considered another extension replacing the base ring $R$ by a semidualizing module $C$. He then introduced and studied the $\GC$-dimension for finitely generated modules over noetherian rings.

Several years later, other generalizations were given in this direction by Holm and J$\phi$rgensen \cite{HJ06}, White \cite{Whi10}, Liu, Huang and Xu  \cite{LHX13} and Xu and Ding \cite{XD16} requiring fewer assumptions, and giving rise to the concept of (not necessarily finitely generated) $\GC$-projective and $\GC$-flat $R$-modules with $C$ being a semidualizing $(R,S)$-bimodule. 

However, requiring $C$ to be a semidualizing module is by no means quite restrictive and some of the conditions that define it have proven, in many cases, to be unnecessary. Therefore, a way to avoid these conditions on $C$ without losing the good properties of the theory has been studied recently in \cite{BGO16a,BEGO22a} where the authors found the minimum conditions to require $C$ to still have a nice theory to develop. Modules satisfying these conditions were called $w$-tilting in the case of $\GC$-projectivity and $w^+$-tilting in the case of $\GC$-flatness.

The purpose of this paper is to continue our investigation on this subject from homological and homotopical aspects. We are mainly interested in the global $\GC$-flat dimension of $R$.

Our first main result is to provide a simple way to compute the global $\GC$-flat dimension of $R$. We prove that the finiteness of this global depends only on the finiteness of the flat dimension of the $\I$-injective right $R$-modules and the $\F$-flat dimension of the injective left $R$-modules. More precisely, we have the following result, which we will prove in Section 3.

\bigskip

\noindent{\bf Theorem A.} \label{thm A} Assume that $_RC$ is $w^+$-tilting and that the class of $\GC$-flat left $R$-modules is closed under extensions. Then, for a positive integer $n$, the following assertions are equivalent:
\begin{enumerate}
	\item $\GCFD(R):=\sup\{\Gfc_R(M)\;|\;M\in R\text{-Mod}\} \leq n.$ 
	\item  The following two assertions hold:
	\begin{enumerate}
		\item  $\fd_R(M)\leq n$ for every  $\I$-injective right $R$-module $M$.
		\item  $\fc_R(M)\leq n$ for every injective left $R$-module $M$.
	\end{enumerate}  
	\item The following two assertions hold:
	\begin{enumerate}
		\item  $\fd_R(M)\leq n$ for every $\I$-injective right $R$-module $M$.
		\item   $\Gfc_R(M)\leq n$ for every injective left $R$-module $M$.    
	\end{enumerate}
\end{enumerate}
$$***$$

It is well known that the weak global dimension of any ring $R$ can be computed by the flat dimension of either the left or right $R$-modules. More precisely,  we have the following equality: 
$$\sup\{\fd_R(M)|M \text{ is a left $R$-module}\}=\sup\{\fd_R(M)|M \text{ is a right $R$-module}\}.$$

The theory of Gorenstein flat dimensions relative to a semidualizing module is usually studied over commutative rings.  But, once we remove the commutativity of the ring, taking into account that the definition of a semidualizing $(R,S)$-bimodule $C$ is left-right symmetric, a question similar to that of the weak global dimension arises: 
\bigskip

\textbf{Question:} Does the following equality hold true?
$$\sup\{\Gfc_R(M)|M\in R\text{-Mod}\}=\sup\{\Gfc_S(M)|M\in \text{Mod-}S\}$$

In Section 4, as an application of Theorem A,  we give a positive answer to this question when the class of $\GC$-flat left $R$-modules and the class of $\GC$-flat right $S$-modules are both closed under extensions (see Theorem \ref{GCF is sym} and Corollary \ref{Sym when R,S coh}).

\bigskip
\noindent{\bf Theorem B.} \label{th B} Let $C$ be a semidualizing $(R,S)$-bimodule and assume that the classes of $\GC$-flat left $R$-modules and $\GC$-flat right $S$-modules are both closed under extensions. Then
$$\sup\{\Gfc_R(M)\;|\;M\in R\text{-Mod}\}=\sup\{\Gfc_S(M)\;|\;M\in \text{Mod-}S\}.$$
In particular, this holds when $R$ is left coherent and $S$ is  right coherent. 

\bigskip

As a direct consequence of Theorem B,  we obtain a positive answer (Corollary \ref{Gfd is sym})  to Bennis' conjecture \cite{Ben10} for any ring $R$.  However, this conjecture was recently solved independently by S. Bouchiba \cite{Bou15} and later by Christensen, Estrada, and Thompson \cite{CET21}. We note that our approach is different from theirs. Theorem B, in particular, provides a new and simpler proof that sheds more light on such symmetries.

$$***$$

In \cite{Hov02} Hovey introduced and studied  abelian  model structures on abelian categories. He showed that an abelian model structure on an abelian category $\A$ is equivalent to a triple of subcategories $(\Q,\W,\R)$ in $\A$ such that $\W$ is thick and $(\Q,\W\cap\R)$ and $(\Q\cap\W,\R)$ are two complete cotorsion pairs. In this case,  $\Q$ (resp., $\W$ and $R$) is the subcategory of $\A$ consisting of all  cofibrant (resp., trivial and fibrant) objects associated to the corresponding abelian model structure. An important feature of an abelian and hereditary model structure is that its homotopy category is  triangulated. Thus, a good and simple way to model triangulated categories is to construct hereditary abelian model structures.  

When $R$ is a commutative noetherian ring and $C$ is a semidualizing module, by a theorem of Hu, Geng, Wu and Li \cite[Theorem 4.3]{HGWL21}, the subcategory $\GCF(R)\cap\B_C(R)\cap\C(R)$ of  $\GC$-flat and $\C$-cotorsion $R$-modules belonging to the Bass class $\B_C(R)$, is a Frobenius category with projective-injective objects all $R$-modules in $\F(R)\cap \C(R)$.  Therefore, the stable category $$\underline{\GCF(R)\cap\B_C(R)\cap\C(R)}$$ is a triangulated category.

From the point of view of homotopy theory, it is natural to ask whether there is a hereditary abelian model structure that models this stable category. In other words, is there a hereditary abelian model structure $\M=(\Q,\W,\R)$ in which its homotopy category $\Ho(\mathcal{M})$ is triangulated equivalent to the stable category $$\underline{\GCF(R)\cap\B_C(R)\cap\C(R)}?$$

The authors of this paper have constructed such a hereditary abelian model structure $\M=(\GCF(R),\W,\B_C(R)\cap\C(R))$ in their recent work \cite{BEGO22b}. This model structure therefore solves the previous problem  but raises a new one: we do not have a clear description of the most important class in any model structure, the class of trivial objects $\W$.

Our next result shows that if the global $\GC$-flat dimension of $R$ is finite, we have a positive answer to this question (see Theorem \ref{GCF model structure}).

\noindent{\bf Theorem C.}\label{Th C} Assume that $C$ is $w$-tilting admitting a degreewise finite projective resolution such that $\GCF(R)$ is closed under extensions. If $R$ has finite global $\GC$-flat dimension, then there exists a hereditary abelian model structure
$$\left(\GCF(R),\widehat{\V_C(R)},\H_C(R)\right),$$
as follows:

$\bullet$ Cofibrations coincide with  monomorphisms with  $\GC$-flat cokernels.

$\bullet$ Fibrations coincide with epimorphisms with $\H_C$-cotorsion kernels.

$\bullet$ Trivial objects coincide with all modules having finite $\V_C$-flat dimension.

\bigskip 

Here, $\H_C(R):=\B_C(R)\cap\C(R)$ and  $\V_C(R)={}^\perp(\B_C(R)\cap\C(R))$ is the left $\Ext$-orthogonal class of $\H_C(R)$.  We call modules in  $\V_C(R)$ and $\H_C(R)$ $\V_C$-flat and $\H_C$-cotorsion modules, respectively. We call them this way as they satisfy most of the nice properties that flat and cotorsion modules have. 

\bigskip


This problem of lack of information about trivial objects arises even in the general setting of abelian categories. In our previous paper \cite{BEGO22b}, in order to construct the $\GC$-flat model structure, we used a well-known result by Gillespie \cite[Theorem 1.1]{Gil15}. This result states that if  $(\Q,\widetilde{\R})$ and $(\widetilde{\Q},\R)$ are two complete hereditary cotorsion pairs in $\A$ such that  $\widetilde{\Q}\subseteq \Q$ and  $\Q\cap \widetilde{\R}=\widetilde{\Q}\cap \R$, then there exists a unique thick class $\W$ such that $(\Q,\W,\R)$ is a Hovey triple. And, the class of  $\W$ is characterized by the classes $\widetilde{\Q}$ and $\widetilde{\R}$. This result is a very useful tool for building model structures and has been used by many researchers for this purpose. 

Researchers interested in model structures would like to know more about trivial objects than what this result provides. The importance of these objects comes mainly from the fact that they determine the associated homotopy category, as explained in the fundamental theorem of model categories \cite[Theorem 1.2.10]{Hov99}.

A new approach is therefore appreciated. For instance, \v{S}aroch and \v{S}\v{t}ov\'{\i}\v{c}ek  \cite[Sections 4 and 5]{SS20}, and later Estrada, Iacob and P\'{e}rez \cite{EIP20} and the second named author \cite{Elm23}, have recently used  new techniques to construct new model structures with an explicit description of the trivial objects.

Our last main result, which is used to prove Theorem C, gives one more step towards a better understanding of this class of trivial objects.
\bigskip

\noindent{\bf Theorem D.} \label{thm D} Assume $(\Q,\widetilde{\R})$ and $(\widetilde{\Q},\R)$ are complete hereditary cotorsion pairs in an abelian category $\A$ such that
\begin{enumerate}
	\item[(a)]  $\widetilde{\Q}\subseteq \Q$ (or equivalently, $\R\subseteq \widetilde{\R}$).
	\item[(b)] $\Q\cap \widetilde{\R}=\widetilde{\Q}\cap \R$.
\end{enumerate}
Then, the following assertions hold:
\begin{enumerate}
	\item 	If $\sup\{\Q\textnormal{-}\resdim_{\A}(A)|A\in \A\}<\infty$, then $(\Q,\W,\R)$ is a Hovey triple where 
	$$\W=\{A\in\A,\; \widetilde{\Q}\textnormal{-}\resdim_{\A}(A)<\infty\}.$$
	
	\item If $\sup\{\R\textnormal{-}\coresdim_{\A}(A)|A\in \A\}<\infty$, then $(\Q,\W,\R)$ is a Hovey triple where 
	$$\W=\{A\in\A,\; \widetilde{\R}\textnormal{-}\coresdim_{\A}(A)<\infty\}.$$
\end{enumerate}

The proof of this result is based on an interesting connection between Hovey triples and (weak) AB-contexts. The later have been studied by Auslander and Buchweitz in \cite{AB89} and named by Hashimoto in his book \cite{Has00}.  See also \cite{BMPS19} for a recent development on this subject. (Weak) AB-contexts are known to be useful to generate special (pre-)covers and (pre-)envelopes. It turns out that they are also useful in constructing Hovey triples (see Theorem \ref{Hov from AB-cont}).

\section{Preliminaries}

\hskip .5cm Throughout this paper, $R$ and $S$ will be  associative (non-necessarily commutative) rings with identity, and all $R$-modules and all $S$-modules will be, unless otherwise specified, unital left $R$-modules and  right $S$-modules, respectively. Right $R$-modules (resp, left $S$-modules) will be identified with left (resp., right) modules over the opposite ring $R^{op}$ (resp, $S^{op}$). We use $\Inj(R)$, $\Proj(R)$, $\Flat(R)$ and $\Cot(R)$ to denote the classes of injective, projective, flat and cotorsion $R$-modules, respectively.

From now on $C$ will stand for an $R$-module, $S$ for its endomorphism ring, $S=\End_R(C)$, and $C^+$  for the character module of $C$, $C^+=\Hom_{\mathbb{Z}}(C,\mathbb{Q}/\mathbb{Z})$. Notice that, $C$ is an $(R,S)$-bimodule and $C^+$ an $(S,R)$-bimodule. We use $\Add_R(C)$ (resp., $\Prod_R(C)$) to denote the class of all $R$-modules which are isomorphic to direct summands of direct sums (resp., direct products) of copies of $C$.

We use $\A$ to denote an abelian category and  by a subcategory of $\A$, we will always mean a full subcategory closed under isomorphisms.  Any class of objects in $\A$ will be thought as a (full) subcategory.

Let $\X$ be a class of objects in $\A$.	A sequence $\XX$ in $\A$ is called $\Hom_R(\X,-)$-exact (resp.,
$ \Hom_R(-,\X)$-exact) if $\Hom_R(X,\XX)$ (resp. $\Hom_R(\XX,X)$)
is an exact complex for every $X\in\X$. In case $\A=R$-Mod and  $\Y$ is a class of right $R$-modules, $\XX$ is called $\left(\Y\otimes_R-\right)$-exact if $Y\otimes_R\XX$
is an exact complex for every $Y\in\Y$. An $\X$-resolution of an object  $A$ is an exact complex $$\cdots\to X_1\to X_0\to A\to 0$$ where
$X_i\in\X$. An $\X$-coresolution of $M$ is defined dually.



The class $\X$  is projectively resolving if it is closed under extensions, kernels of epimorphisms, and contains the projective objects of $\A$.  $\X$ is left thick if it is closed under extensions, kernels of epimorphisms, and direct summands. Injectively coresolving and right thick classes are defined dually. $\X$ is thick if it is left and right thick.
\bigskip

\noindent{\bf (Co)resolution dimension.} An object $A\in\A$ is said to have $\X$-resolution dimension less than or equal to an integer $n\geq 0$, $\X\textnormal{-}\resdim_\A(A)\leq n$, if $A$ has a finite $\X$-resolution: $0\to X_n\to\cdots\to X_1\to X_0\to M\to 0.$

If $n$ is the least non negative integer for which such a sequence exists then its $\X$-resolution dimension is precisely $n,$ and if there is no such $n$ then we define its $\X$-resolution dimension as $\infty$. We denote by $\widehat{\X}$ the full subcategory of objects in $\A$ having a finite $\X$-resolution. The $\X$-resolution dimension of a class $\Y\subseteq \A$, is defined as 
$$\X\textnormal{-}\resdim_\A(\Y)=\sup\{\X\textnormal{-}\resdim_\A(Y)| Y\in \Y\}.$$

\noindent{\bf Cotorsion pairs.} Given a class $\X\subseteq \A$, we use  $\X^{\perp}$ to denote the class of all objects $M\in \A$ such that $\Ext^1_\A(X,M)=0$  for all $X\in\X$. Similarly,  $^{\perp}\X=\{M\in \A|\Ext^1_\A(M,X)=0, \forall X\in\X\}$.



A pair $\left(\X,\Y\right)$ of classes of objects in $\A$ is called a cotorsion pair if $\X^\perp=\Y$ and $\X={}^\perp\Y$. A cotorsion pair $\left(\X,\Y\right)$ is said to be hereditary if $\Ext^i_\A(\X,\Y)=0$ for every $i\geq 1$. 
A cotorsion pair $(\X,\Y)$ is called complete if, for any $A\in\A$, there exist two short exact sequences
$$0\to Y\to X\to A\to 0 \text{ and } 0\to A\to Y'\to X'\to 0$$
with $X,X'\in \X$ and $Y,Y'\in\Y$.


\noindent{\bf Abelian model structures.} Assume that $\A$ is a bicomplete category. A model structure
on $\A$ is given by three classes of morphisms of $\A$, called cofibrations, fibrations and weak equivalences, that satisfy a set of axioms that enable one to introduce a homotopy theory on that category (\cite[Definition 1.1.3.]{Hov99}).

Hovey defined in \cite{Hov02} the notion of an abelian model structure and showed \cite[Theorem 2.2]{Hov02} that there is a close link between abelian model structures and some cotorsion pairs. Namely, an abelian model structure on any abelian category $\A$ is equivalent to a triple $\mathcal{M}=(\Q,\W,\R)$ (known as a Hovey triple) of classes of objects in $\A$ such that $\mathcal{W}$ is thick and $(\mathcal{Q},\mathcal{W}\cap\mathcal{R})$ and $(\mathcal{Q}\cap\mathcal{W 
},\mathcal{R})$ are complete cotorsion pairs.  In this case, $\Q$ is the class of cofibrant objects, $\R$ is the class of fibrant objects, and $\W$ is the class of trivial objects of the model structure. Here, an abelian model structure (\cite{Hov02}, see also \cite{Gil16}) is a model structure  over a bicomplete abelian category $\A$ such that:  (1) a morphism is a (trivial) cofibration if and only if it is a monomorphism with
(trivially) cofibrant cokernel and (2) a morphism is a (trivial) fibration if and only if it is an epimorphism with (trivially) fibrant kernel.

A Hovey triple is hereditary if both of the associated cotorsion pairs are hereditary. In this paper, we often identify any abelian model structure with its associated Hovey triple.
\bigskip

\noindent{\bf (Weak) Auslander-Buchweitz contexts.} 
Given two subcategories $\omega$ and  $\X$ of $\A$, $\omega$ is cogenerator for $\X$ if  $\omega \subseteq \X$ and for any $X\in\X$ there exists an exact sequence $0\to X\to W\to X'\to 0$, such that $X'\in\X$ and $W\in\omega$. If moreover, $\Ext^i_\A(\X,\omega)=0$ for all $i\geq 1$, $\omega$ is said to be an injective cogenerator for $\X$. 

Recall from \cite{Has00} that a triplet $(\X,\Y,\omega)$ of subcategories of $\A$ is a left weak Auslander-Buchweitz context (or a weak AB-context for short), in the terminology of Hashimoto, if the following three conditions are satisfied:

\begin{enumerate}
	\item[(AB1)]  $\X$ is left thick.
	
	\item[(AB2)] $\Y$ is right thick and $\mathcal{Y}\subseteq \widehat{\mathcal{X}}$.
	
	\item[(AB3)] $\omega=\X\cap\Y$, and $\omega$ is an injective cogenerator for $\X$.
\end{enumerate}
If moreover $\widehat{\X}=\A$, then it is called a left Auslander-Buchweitz context (or AB-context for short). 

\bigskip

\noindent \textbf{Relative Gorenstein flat dimensions.} 
\begin{defn} A right $R$-module $N$ is said to be $\I$-injective if it is in $\Prod_{R^{op}}(C^+)$ and an $R$-module $M$ is said to be $\F$-flat if its character module, $M^+$, is $\I$-injective. 
	
	We let $\F(R)$ and  $\I(R^{op})$ denote the classes of left $\F$-flat and right $\I$-injective modules, respectively.
\end{defn}

$\F$-flat and $\I$-injective modules are natural generalizations of flat and injective modules to the relative setting. They also generalize $C$-flat  and $C$-injective modules in the sense of  \cite[Definition 5.1]{HW06}.
\begin{exmps}\label{exmps of Fc-flat} The following examples are used repeatedly throughout the paper.
	\begin{enumerate}
		\item If $_RC=R$ (or more generally, $_RC$ is a flat generator of $R$-Mod), then 
		$$\F(R)=\Flat(R) \text{ and } \I(R^{op})=\Inj(R^{op}).$$
		\item \cite[Proposition 3.3]{BEGO22a} If $_RC$ is finitely presented, then
		$$\F(R)=C\otimes_S\Flat(S^{op}) \text{ and  } \I(R^{op})=\Hom_S(C,\Inj(S)).$$
		
	\end{enumerate}
\end{exmps}


\begin{defn} [\cite{BEGO22a}, Definition 4.1] An $R$-module $C$ is said to be $w^+$-tilting if it satisfies the following two properties:
	\begin{enumerate}
		\item $C$ is $\prod$-$\Tor$-orthogonal, that is, $\Tor^R_{i\geq 1}((C^+)^I,C)=0$ for every set $I$.
		\item There exists a  $\left( \I(R^{op})\otimes_R-\right) $-exact $\F(R)$-coresolution $$\textbf{X}:0\rightarrow R\rightarrow C_{-1}\rightarrow C_{-2}\rightarrow\cdots$$
	\end{enumerate}
\end{defn}

\begin{defn} [\cite{BEGO22a}, Definition 4.4]\label{GcF-def} An $R$-module $M$ is said to be $\GC$-flat if there exists an exact and $\left(\I(R^{op})\otimes_R-\right)$-exact sequence $$\textbf{X}:\cdots \rightarrow F_1\rightarrow F_0 \rightarrow C_{-1}\rightarrow C_{-2}\rightarrow\cdots$$ with each $C_i\in\F(R)$ and $F_j\in\Flat(R)$, such that $M\cong \Im(F_0\to C_{-1})$.
	
	The class of all $\GC$-flat $R$-modules is denoted by $\GCF(R)$. 
\end{defn}	
\begin{defn}[\cite{BEGO22b}] An $R$-module $M$ is said to be $\GC$-cotorsion (resp., $\C$-cotorsion) if $\Ext_R^1(N,M)=0$ for all $\GC$-flat (resp., $\F$-flat) $R$-modules $N$.
	
	The class of all $\GC$-cotorsion (resp., $\C$-cotorsion) $R$-modules is denoted by $\GCC(R)$ (resp., $\C(R)$). 
\end{defn}
\begin{defn}[\cite{BEGO22a}, Definition 4.4] The $\F$-flat and $\GC$-flat dimension of an $R$-module $M$ are defined, respectively, as: 
	$$\fc_R(M)=\F(R)\text{-}\resdim_R(M)\text{ and }\Gfc_R(M)=\GCF(R)\text{-}\resdim_R(M).$$
\end{defn}

Some of our results are related to the condition that the class of $\GC$-flat modules is closed under extensions. 
\begin{defn}[\cite{BEGO22a}, Definition 4.14]  A ring $R$ is said to be $\GCF$-closed provided that the class of $\GC$-flat $R$-modules is closed under extensions.
\end{defn} 

The following  lemma provides different situations in which the class of $\GC$-flat $R$-modules is closed under extensions, which are taken from \cite{BEGO22a} and \cite{BEGO22b}.
\begin{lem}\label{GCF-closed situations} Assume that $C$ is $\prod$-$\Tor$-orthogonal. The ring $R$ is $\GCF$-closed in the following cases:
	\begin{enumerate}
		\item[(a)] $_RC$ is a flat generator.
		\item[(b)] $\F(R)$ is closed under direct products. For instance, if $S=\End_R(C)$ is right coherent and both $_RC$ and $C_S$ are finitely presented.
		\item[(c)] Every $R$-module has finite $\F$-flat dimension.
	\end{enumerate}
\end{lem}
\proof (a)  Since $_RC$ is a flat generator, $\F(R)=\Flat(R)$. It is immediate to check that in \cite[Proposition 4.23(2)$\Rightarrow$ (1)]{BEGO22a} the assumptions that $R$ is ${\rm G_WF}$-closed and  $W$ is $w^+$-tilting are not necessary.  The assertion, therefore, follows by the proof of \cite[Proposition 4.23(2)$\Rightarrow$ (1)]{BEGO22a}  and \cite[Theorem 4.11(3)]{SS20}.

(b) Follows by \cite[Corollary 4.13]{BEGO22a} and \cite[Lemma 4.7]{BEGO22b}.

(c) Follows by \cite[Corollary 7.5]{BEGO22a}.
\cqfd

Finally, we recall the definition of Foxby classes. The Bass class will be of interest at the end of this paper.

Associated to the bimodule $_RC_S$, we have the Auslander and Bass classes, $\A_C(S)\subseteq S$-Mod and $\B_C(R)\subseteq R$-Mod, respectively, defined as follows:

$\bullet$ $\A_C(S)$ is the class of all left $S$-modules $M$ satisfying:

$\Tor^S_{\geq 1}(C,M)=0=\Ext_R^{\geq 1}(C,C\otimes_S M)$ and the canonical map $$\mu_M:M\to \Hom_{S^{op}}(C,C\otimes_RM)$$ is an isomorphism of left $S$-modules.

$\bullet$ $\B_C(R)$ consists of all $R$-modules $N$ satisfying:

$\Ext_R^{\geq 1}(C,N)=0=\Tor^S_{\geq 1}(C,\Hom_R(C,N))=0$ and the canonical map $$\nu_N:C\otimes_S\Hom_R(C,N)\to N$$ is an isomorphism of $R$-modules.

\section{Relative weak Gorenstein global dimension} 

\hskip .5cm In this section we define and study the global $\GC$-flat dimension of $R$.

\begin{defn} The global $\GC$-flat  dimension of $R$ is defined as the supremum, if it exists, of the $\GC$-flat dimension of every $R$-module:
	$$\GCFD(R):=\sup\{\Gfc_R(M)|\text{ M is an $R$-module }  \}.$$
	We set $\GCFD(R)=\infty$ if such a supremum does not exist. 
\end{defn}

\begin{rem}\begin{enumerate}
		\item The global $\GC$-flat dimension of rings was briefly studied in \cite{ZS13,ZW18} over commutative rings with $C$ being a semidualizing $R$-module. 
		\item   When $_RC$ is a flat generator, we recover $\GCFD(R)$ as $\GFD(R)$, the weak global Gorenstein  dimension of $R$.
	\end{enumerate}
\end{rem}

In order to prove Theorem A, we will need the following two lemmas.

\begin{lem}\label{nec-cond of when Gc-fd is finite} Let $R$ be $\GCF$-closed and $_RC$ be  $\prod$-$\Tor$-orthogonal. For an $R$-module $M$ and an integer $n\geq 1$, if $\Gfc_R(M)\leq n,$ then there exists an exact sequence of $R$-modules $0\rightarrow F\rightarrow G\rightarrow M\rightarrow 0$, where $G$ is $\GC$-flat and $\fc_R(F)\leq n-1.$
	
\end{lem}
\proof By definition, there is a short exact sequence of $R$-modules $0\to K\to G_0\to M\to 0$, where $G_0$ is $\GC$-flat  and $\Gfc_R(K)\leq n-1$. By \cite[Proposition 7.11]{BEGO22a}, there is an exact sequence of $R$-modules $0\rightarrow K\rightarrow F\rightarrow L\rightarrow 0$, where $L$ is $\GC$-flat and $\fc_R(F)\leq n-1.$ 

Consider the following pushout diagram:

$$\xymatrix{ & 0\ar[d] & 0\ar[d] & & \\  0\ar[r] & K\ar[d] \ar[r] & G_0\ar@{-->}[d] \ar[r] & M \ar@{=}[d] \ar[r] & 0 \\ 0\ar[r] & F \ar@{-->}[r]\ar[d] & G \ar[d] \ar[r] & M\ar[r] & 0 \\ & L\ar@{=}[r] \ar[d] & L\ar[d] & & \\ & 0 & 0}.$$
Therefore, $G$ is $\GC$-flat as $G_0$ and $L$ are $\GC$-flat and $R$ is $\GCF$-closed.
\cqfd

In the absolute case \cite[Corollary 2.3]{Ben10}, the following key lemma is based on the fact that the class of modules of finite flat dimension is closed under direct summands which is not the case in our relative setting. Here, we adopt a different proof.

\begin{lem}\label{Gc-fd of inj} Let $R$ be $\GCF$-closed and $C$ be a $w^+$-tilting $R$-module. If $M$ is an injective $R$-module, then $\fc_R(M)=\Gfc_R(M)$.
\end{lem}
\proof The inequality $\Gfc_R(M) \leq \fc_R(M)$ holds by \cite[Theorem 7.12]{BEGO22a}. For the other inequality, we may assume that $n=\Gfc_R(M)<\infty$. 

If $n=0$ then $M$ is $\GC$-flat and since $M$ is injective there exists by \cite[Proposition 4.10]{BEGO22a} a split exact sequence $0\to M\to V\to L\to 0$ with $V\in \F(R)$. Thus, $M\in\F(R)$ and hence $\fc_R(M)=0$. 

Let us now assume that  $n\geq 1$. By \cite[Proposition 7.11]{BEGO22a} and Lemma \ref{nec-cond of when Gc-fd is finite}, there exist two exact sequences of $R$-modules 
$$0\to M\to F_1\to G_1\to 0 \text{ and }  0\to F_2\to G_2\to M\to 0$$
where $G_1$ and $G_2$ are $\GC$-flat and $\fc_R(F_1)\leq n$ and $\fc_R(F_2)\leq n-1$. Since $M$ is injective, The first sequence splits and so $M\oplus G_1\cong F_1$. Then, adding the second sequence with $0\to 0\to G_1\to G_1\to 0$ we get a short exact sequence of the form  $0\to F_2\to G_2\oplus G_1\to F_1\to 0$. 

Since $\fc_R(F_1)\leq n$ and $\fc_R(F_2)\leq n-1$, there exist exact sequences 
$$0\to X_{n-1}\to\cdots\to X_0\to F_1\to 0\text{ and }0\to Y_n\to Y_{n-1}\to\cdots\to Y_0\to F_2\to 0$$
where $X_i,Y_j\in\F(R)$. It follows by \cite[Proposition 3.5]{BEGO22a} that $C^+$ is $\prod$-self-orthogonal. This implies that $\Ext_{R^{op}}^{k\geq 1}(Y_j^+,E)=0$ for all $j=0,\cdots,n$ and all $E\in \I(R^{op})$, so by dimension shifting we get that $\Ext_{R^{op}}^1(F_2^+,E)=0$. Hence, the exact sequence $$ 0 \to F_1^+ \to  G_1^+\oplus G_2^+ \to F_2^+\to 0$$ is $\Hom_{R}(-,\I(R^{op}))$-exact and by Horseshoe Lemma we get the commutative diagram with exact rows and columns:

$$\xymatrix{ & 0 \ar[d] & 0 \ar[d] & 0 \ar[d] \\ 0 \ar[r] & F_1^+ \ar[r] \ar[d] & G_1^+\oplus G_2^+ \ar[r] \ar[d] & F_2^+ \ar[r] \ar[d] & 0 \\ 0 \ar[r] & X_0^+ \ar[r] \ar[d] & X_0^+\oplus Y_0^+ \ar[r] \ar[d] & Y_0^+ \ar[r] \ar[d] & 0 \\ & \vdots \ar[d] & \vdots \ar[d] & \vdots \ar[d] \\  0 \ar[r] & X_{n-1}^+ \ar[r] \ar[d] & X_{n-1}^+\oplus Y_{n-1}^+ \ar[r] \ar[d] & Y_{n-1}^+ \ar[r] \ar[d] & 0 \\ 0 \ar[r] & 0 \ar[r] \ar[d] &  K \ar[r] \ar[d] & Y_n^+ \ar[r] \ar[d] & 0 \\ & 0 & 0 & 0}.$$
From the diagram above, one can see that $K\cong Y_n^+\in\I(R^{op})$. Then, we get $\ic_{R^{op}}(G_1^+\oplus G_2^+)\leq n$. By \cite[Proposition 4.7 and Theorem 4.12]{BEGO22a}, $(G_1\oplus G_2)^+$ is $\rm G_{C^+}$-injective, since $G_1$ and $G_2$ are $\GC$-flat. This implies by \cite[Proposition 3.4(2)]{BGO16b}, that  $(G_1\oplus G_2)^+ \cong G_1^+\oplus G_2^+\in \I(R^{op})$ which gives that $G_2\in\F(R)$. Thus, $\fc_R(M)\leq n=\Gfc_R(M)$. \cqfd

Now we are ready to prove  Theorem A.	

\begin{thm} \label{charc of glob Gc-flat dim} Assume that $R$ is $\GCF$-closed and $_RC$ is $w^+$-tilting. Then, for a positive integer $n$, the following assertions are equivalent:
	\begin{enumerate}
		\item $\GCFD(R)\leq n.$ 
		\item  The following two assertions hold:
		\begin{enumerate}
			\item  $\fd_{R^{op}}(M)\leq n$ for every  $\I$-injective right $R$-module $M$.
			\item  $\fc_R(M)\leq n$ for every injective left $R$-module $M$.
		\end{enumerate}  
		\item The following two assertions hold:
		\begin{enumerate}
			\item  $\fd_{R^{op}}(M)\leq n$ for every $\I$-injective right $R$-module $M$.
			\item   $\Gfc_R(M)\leq n$ for every injective left $R$-module $M$.    
		\end{enumerate}
	\end{enumerate}
	
	Consequently, the global $\GC$-flat  dimension of $R$  can be computed via the following formula: 
	$$\GCFD(R)=max\{\fd_{R^{op}}(\I(R^{op})),\fc_R(\Inj(R))\}.$$
\end{thm} 
\proof $2.\Leftrightarrow 3.$ Follows by Lemma  \ref{Gc-fd of inj}.

$1.\Rightarrow 3.$ The assertion $(b)$ holds by definition. Let us prove $(a)$. 

Note that $\fd_{R^{op}}(M)\leq n$ if and only if $\Tor^R_{n+1}(M,A)=0$ for every $R$-module $A$. But $\GCFD(R)\leq n$ means that for every such $A$ there is an exact sequence
$$0\to G_n\to\cdots\to G_0\to A\to 0$$
with $G_i\in \GCF(R),$ for all $i\geq 0$. Now, $M\in \I(R^{op})$ implies that $\Tor^R_{i\geq 1}(M,G_i)=0,$ for all $i$ so $\Tor^R_{n+1}(M,A)\cong \Tor^R_1(M,G_n)=0$.

$3.\Rightarrow 1.$ Let $M$ be an $R$-module. Consider a projective resolution and an injective coresolution of $M$:$$\cdots \to P_1\to P_0\to P_{-1}=M\to 0\text{ \;\;and\;\;  }0\to M=I_{-1}\to I_0\to I_1\to\cdots,$$ 
respectively. Decomposing these exact sequences into short exact ones we get, for every integer $i\in \mathbb{N}$, 
$$0 \to N_{i+1}\to P_i\to N_i\to 0\text{ \;\;and\;\;  }0\to K_i\to I_i\to K_{i+1}\to 0$$ 
where $N_i=\Im(P_i\to P_{i-1})$ and  $K_i=\Ker(I_i\to I_{i+1})$.  
Adding the direct sum of the first sequences, 
$$0 \to \oplus_{i\in\mathbb{N}}N_{i+1}\to \oplus_{i\in\mathbb{N}}P_i\to M\oplus( \oplus_{i\in\mathbb{N}}N_{i+1})\to 0,$$
to the direct product of the second ones, 
$$0\to  M\oplus( \prod_{i\in\mathbb{N}}K_{i+1})\to \prod_{i\in\mathbb{N}}I_i\to \prod_{i\in\mathbb{N}}K_{i+1}\to 0,$$
we get the  exact sequence 
$$0\to M\oplus (N\oplus K)\to P\oplus I\to M\oplus(N\oplus K)\to 0$$
where $N=\oplus_{i\in \mathbb{N}} N_{i+1}$, $K=\prod_{i\in \mathbb{N}} K_{i+1}$, $P=\oplus_{i\in \mathbb{N}} P_i$ and $I= \prod_{i\in \mathbb{N}} I_i$.

Now, consider a projective resolution of $M\oplus (N\oplus K):$

$$\cdots\to H_1 \to H_0\to M\oplus (N\oplus K)\to 0.$$
Thus, by Horseshoe Lemma, we get a commutative diagram with exact rows and columns
$$\xymatrix{ 	&  &0\ar[d]  &0\ar[d]   &0\ar[d] \\
	&0\ar[r]  & J_n\ar[r]\ar[d] & H\ar[r]\ar[d] & J_n\ar[r]\ar[d] &0\\
	&0\ar[r]  & H_{n-1}\ar[r]\ar[d] & H_{n-1}\oplus H_{n-1}\ar[r]\ar[d] & H_{n-1}\ar[r]\ar[d] &0\\
	&  & \colon \ar[d] & \colon\ar[d] & \colon\ar[d] \\
	&0\ar[r]  & H_0\ar[r]\ar[d] & H_0\oplus H_0\ar[r]\ar[d] & H_0\ar[r]\ar[d] &0\\
	&0\ar[r]  & M\oplus (N\oplus K)\ar[r]\ar[d] & P\oplus I\ar[r]\ar[d] & M\oplus(N\oplus K)\ar[r]\ar[d] &0\\
	&  &0  &0   &0
}$$

Note that $P$ is $\GC$-flat by \cite[Proposition 4.17]{BEGO22a} and since $I$ is injective,  $\Gfc_R(P\oplus I)=\Gfc_R(I)\leq n$ by \cite[Corollary 7.8(2)]{BEGO22a} and the hypotheses.  But, since each $H_i$ is $\GC$-flat, $H$ is $\GC$-flat as well by \cite[Theorem 7.7]{BEGO22a}. 

Now, let $X_R$ be any $\I$-injective right $R$-module. By hypothesis, $\fd_{R^{op}}(X)\leq n$. Then, using the projective resolution of $M\oplus (N\oplus K)$, we get 
$$\Tor_1^R(X,J_n)\cong \Tor_{n+1}^R(X,M\oplus (N\oplus K))=0.$$
So, the sequence $0\to J_n\to H\to J_n\to 0$ is $\left(\I(R^{op})\otimes_R-\right)$-exact. Then, assembling these sequences, we get a $\left(\I(R^{op})\otimes_R-\right)$-exact exact sequence  
$$\cdots\to H\stackrel{f}{\to} H\stackrel{f}{\to} H\to \cdots$$ with $J_n=\Ker f$. Therefore,  $J_n$ is $\GC$-flat by \cite[Theorem 5.1]{BEGO22a} and then $\Gfc_R(M)\leq \Gfc_R(M\oplus (N\oplus K))\leq n$ by \cite[Corollary 7.8]{BEGO22a} as desired. \cqfd

The following special case of Theorem \ref{charc of glob Gc-flat dim} was proved by Emmanouil (\cite[Theorem 5.3]{Emm12}) when the Gorenstein weak dimension of $R$ is finite. Here we drop this finiteness condition.
\begin{cor}  The weak global Gorenstein  dimension of $R$ can be computed via the  following simple formulas:
	\begin{eqnarray*}
		\GFD(R) &=& max\{\fd_{R^{op}}(\Inj(R^{op})),\Gfd_R(\Inj(R))\}\\
		&=& max\{\fd_{R^{op}}(\Inj(R^{op})),\fd_R(\Inj(R))\}.
	\end{eqnarray*}
\end{cor}

\section{Left-right symmetry of global $\GC$-flat dimension}
\hskip .5cm This section is devoted to prove Theorem B and some of its consequences.

\begin{defn}[\cite{HW06}, Definition 2.1] An $(R,S)$-bimodule $C$ is semidualizing if:
	\begin{enumerate}
		\item $_RC$ and $C_S$ both admit a degreewise finite projective resolution.
		\item $\Ext_R^{\geq 1}(C,C)=\Ext_S^{\geq 1}(C,C)=0.$
		\item The natural homothety maps $R\rightarrow\Hom_S(C,C)$ and $S \rightarrow \Hom_R(C,C)$  are both ring isomorphisms.
	\end{enumerate}
\end{defn}

Here by a degreewise finite projective resolution  we mean a projective resolution in which  every projective module is finitely generated.

For the rest of this section, $C$ will be, unless otherwise stated, a semidualizing $(R,S)$-bimodule.

\begin{lem} \label{flat and Fc-flat dim} The following equalities hold:
	\begin{enumerate}
		\item $\fd_{S^{op}}(\Hom_R(C,E))=\fc_R(E)$ for every injective $R$-module.
		
		Consequently, $\fd_{S^{op}}(\I(S^{op}))=\fc_R(\Inj(R))$.
		\item $\fd_R(\Hom_S(C,E))=\fc_S(E)$ for every injective $S$-module.
		
		Consequently, $\fd_{R^{op}}(\I(R^{op}))=\fc_{S}(\Inj(S))$.
	\end{enumerate}
\end{lem}
\proof  Assertion (2) can be proved in a similar way to (1); so we only prove (1).

First, we prove that  $\fd_{S^{op}}(\Hom_R(C,E))\leq \fc_R(E)$. 

We may assume that $n=\fc_R(E)<\infty$ since the infinite case is clear. Recall (Example \ref{exmps of Fc-flat}(2)) that with our hypothesis we have $\F(R)=C\otimes_S \Flat(S^{op})$, so there exists an exact sequence of $R$-modules:
$$0\to C\otimes_SF_n\to\cdots\to C\otimes_SF_0\to E\to 0$$
where each $_SF_i$ is flat. By \cite[Theorem 3.2.15 and Remark 3.2.27]{EJ00}, for every $k\geq 1$ we get
\begin{eqnarray*}
	\Ext_R^k(C,C\otimes_SF_i)^+&\cong& (\Ext_R^k(C,C)\otimes_S F_i)^+\\
	&\cong& \Hom_{S^{op}}(F,\Ext_R^k(C,C)^+)=0
\end{eqnarray*} 
Then, $\Ext_R^k(C,C\otimes_SF_i)=0$ for every $i\geq 1$. Hence, the sequence of left $S$-modules 
$$0\to \Hom_R(C,C\otimes_SF_n)\to\cdots\to \Hom_R(C,C\otimes_SF_0)\to \Hom_R(C,E)\to 0$$
is exact. But, for each $i=0,\cdots, n$ $$\Hom_R(C,C\otimes_SF_i)\cong \Hom_R(C,C)\otimes_S F_i\cong \;_SF_i$$  by \cite[Theorem 3.2.14]{EJ00}, which implies that $\fd_{S^{op}}(\Hom_R(C,E))\leq n.$

Now we prove the other inequality. Suppose that $n=\fd_{S^{op}}(\Hom_R(C,E)) <\infty$. Then, there exists a finite flat resolution of $_S\Hom_R(C,E)$: 
$$0\to F_n\to \cdots \to F_0\to \Hom_R(C,E)\to 0.$$
Using \cite[Lemma 1.2.11]{GT12}, we get that, for every $k\geq 0$, $$\Tor^S_k(C,\Hom_R(C,E))\cong \Hom_R(\Ext^k_S(C,C),E)\cong \begin{cases}E \text{ if } k=0\\[0.5cm]
0 \text{ if } k\geq 1 \end{cases}.$$ Then, we get an exact sequence
$$0\to C\otimes_SF_n\to \cdots \to C\otimes_S F_0\to E\to 0$$
where each $C\otimes_S F_i$ is $\F$-flat. Hence, $\fc_R(E)\leq n=\fd_{S^{op}}(\Hom_R(C,E)).$ 

Finally, keeping in mind that $\I(S^{op})=\Hom_R(C,\Inj(R))$, we get the last equality. \cqfd

For completeness, we state the right version of Theorem \ref{charc of glob Gc-flat dim}.
\begin{thm}\label {Thm A for S} Assume that $S$ is any $\GCF$-closed ring and $C$ is a $w^+$-tilting $S$-module.  The global $\GC$-flat  dimension of $S$ can be computed as follows: 
	\begin{eqnarray*}
		\GCFD(S)	&=&   max\{\fd_{{S^{op}}}(\I(S^{op})),\fc_S(\Inj(S))\}\\ &=& max\{\fd_{S^{op}}(\I(S^{op})),\Gfc_S(\Inj(S))\}.
	\end{eqnarray*}
	
\end{thm} 
Now, we are in a position to prove Theorem B of the Introduction.
\begin{thm} \label{GCF is sym}Let $R$ be left $\GCF$-closed and $S$ be right $\GCF$-closed. Then
	
	$$\GCFD(R)=\GCFD(S),$$
	
	In this case, we define the common value of these two	numbers to be the $\GC$-flat dimension of the pair of rings $(R,S)$ and denote it by $\GCFD(R,S).$ 
\end{thm}
\proof  By Theorems \ref{charc of glob Gc-flat dim} and \ref{Thm A for S}  we have the following two formulas: $$\GCFD(R)= max\{\fd_{R^{op}}(\I(R^{op})),\fc_R(\Inj(R))\},$$
$$\GCFD(S)=max\{\fd_{S^{op}}(\I(S^{op})),\fc_S(\Inj(S))\}.$$
On the other hand, by Lemma \ref{flat and Fc-flat dim} we know that $$\fd_{R^{op}}(\I(R^{op}))=\fc_S(\Inj(S))\text{ and }  \fd_{R^{op}}(\I(S^{op}))=\fc_R(\Inj(R)).$$
Therefore, $\GCFD(R)=\GCFD(S)$.
\cqfd

We end this section with some consequences of Theorem \ref{GCF is sym}.
\bigskip


Substituting the $\GCF$-closeness condition in  Theorem \ref{GCF is sym} by the coherence condition, keeping in mind Lemma \ref{GCF-closed situations}(b), we get the following special case.
\begin{cor}\label{Sym when R,S coh}Let $R$ be left coherent and $S$ be right coherent. Then
	$$\GCFD(R)=\GCFD(S).$$
	
\end{cor}
\bigskip

Bennis conjectured in \cite{Ben10} that the weak global Gorenstein  dimension of $R$ is symmetric. That is, the equality $\GFD(R)=\GFD(R^{op})$ holds. 

This conjecture has been investigated by many authors. In \cite[Corollary 2.8]{MT09}, Mahdou and Tamekkante  solved it in the case when $R$ is (two-sided) coherent and Emmanouil \cite[Theorem 5.3]{Emm12} solved it in the case when both $R$ and $R^{op}$ have finite Gorenstein weak global dimension. Bouchiba \cite[Corollary 3.5]{Bou15}, on the other hand, showed that this conjecture is true when the classes of Gorenstein flat left and right $R$-modules are closed under extensions. But, this is always the case by the work of \v{S}aroch and  \v{S}\v{t}ov\'{\i}\v{c}ek \cite{SS20} (see Lemma \ref{GCF-closed situations}(a)). Thus, this conjecture is solved. Recently,  Christensen, Estrada, and Thompson have re-established this fact using a different approach \cite[Corollary 2.5]{CET21}. 

As a direct consequence of Theorem \ref{GCF is sym}, we obtain a third proof.  Note that our proof is simple and completely different from the ones given in \cite{Bou15} and \cite{CET21}.
\begin{cor}(\cite[Corollary 3.5.]{Bou15} and \cite[Corollary 2.5]{CET21})\label{Gfd is sym}  Over any ring $R$, the weak global Gorenstein  dimension of $R$ is symmetric, that is, we have the equality:
	$$\GFD(R)=\GFD(R^{op}).$$
	
\end{cor}

\bigskip

Recall \cite{Ben10} that  a ring $R$ is $n$-IF if every left and right injective  $R$-module has flat dimension at most $n$. A two-sided noetherian ring is  $n$-IF if and only if it is $n$-Gorenstein  by \cite[Theorem 9.1.11]{EJ00}. Bennis characterized \cite[Theorem 2.8]{Ben10} $n$-IF rings provided that they are (two-sided) coherent. As another consequence of Theorem \ref{charc of glob Gc-flat dim},  we drop the coherence assumption. 
\begin{cor} A ring $R$ is $n$-IF if and only if $\GFD(R)\leq n$.
\end{cor}

Once the symmetry of the Gorenstein weak global dimension has been established, this corollary can also be deduced from \cite[Theorem 5.3]{Emm12}. 

\bigskip

Assume that  $R$ is left coherent and $S$ is right coherent. Zhu and Ding (\cite[Theorem 2.6]{ZD09}) proved that if $_RC$ and $C_S$ have finite FP-injective dimension, then $\FPid_R(C)=\FPid_S(C).$ We show next that this happens exactly when $R$ (or $S$) has finite global $\GC$-flat dimension. But first, we need the following lemma.
\begin{lem}\label{flat dim and F-inj dim}\item
	\begin{enumerate}
		\item If $R$ is left coherent, then  $\fd_{R^{op}}(\I(R^{op}))=\FPid_R(C)$.
		\item If $S$ is right coherent, then $\fd_{S^{op}}(\I(S^{op}))=\FPid_S(C)$.
	\end{enumerate}
\end{lem}
\proof We only prove (1) since (2) has a similar proof.  

Follwoing \cite[Theorem 2.2]{Fie72}, we get that  $\FPid_R(C)=\fd_{R^{op}}(C^+)\leq \fd_{R^{op}}(\I(R^{op}))$. Conversely, if $\FPid_R(C)=\infty$, then equality holds true. So, we may assume that $n=\FPid_R(C)<\infty$. 

Let $X\in\I(R^{op})$. Then, $X$ is a direct summand of some $(C^+)^I$. Using again \cite[Theorem 2.2]{Fie72}, we get that  $\fd_{R^{op}}(C^+)=\FPid_R(C)=n$. But, the direct product of any flat resolution of a right $R$-module is a flat resolution by \cite[Theorem 3.2.24]{EJ00}. This means that $\fd_{R^{op}}((C^+)^I)\leq n$. Hence,  $\fd_{R^{op}}(X)\leq \fd_{R^{op}}((C^+)^I)\leq n$. Consequently,  $\fd_{R^{op}}(\I(R^{op}))\leq n$. \cqfd

Assume that $R$ and $S$ are left and right noetherian rings, respectively. Recall ( \cite[Definition 3.1]{EJL05}) that a semidualizing $(R,S)$-bimodule is called dualizing if $C_S$ and $_RC$ both have finite injective dimension. Replacing noetherian by coherent and injective by FP-injective we get the following weaker notion.

\begin{defn}Let $R$ and $S$ be left and right coherent rings, respectively. A semidualizing $(R,S)$-bimodule $C$ is said to be a weak dualizing module if $_RC$ and $C_S$ both have finite FP-injective dimension. In this case, $C$ is said to be weak $n$-dualizing where $n=\FPid_R(C)=\FPid_S(C)$.
\end{defn}
\begin{cor} Assume that $R$ is left coherent and $S$ is right coherent. Then,  $\GCFD(R,S)\leq n$ if and only if  $C$ is weak $n$-dualizing.
	
	Consequently, $\GCFD(R,S)=max\{\FPid_R(C),\FPid_S(C)\}.$
\end{cor}
\proof Follows from Corollary \ref{Sym when R,S coh}, Theorem \ref{charc of glob Gc-flat dim},  Lemmas \ref{flat and Fc-flat dim} and \ref{flat dim and F-inj dim} and  the comment just before Lemma \ref{flat dim and F-inj dim}.\cqfd
\begin{cor} Assume that $R$ is left coherent and $S$ is right coherent. Then,  $\GCFD(R,S)=0$ if and only if both $_RC$ and $C_S$ are FP-injective modules.
\end{cor}

\section{$\GC$-flat model structure}

\hskip .5cm In this section we describe the class of weak equivalences in the $\GC$-flat model structure recently found in \cite[Section 4]{BEGO22b}, under the finiteness of the global $\GC$-flat dimension of $R$. This result is a consequence of a more general result (Theorem D of the Introduction) that we will prove in this section. First, we need to develop some preliminary results that we will use throughout this section.
\bigskip

Recall that a class $\X$ of objects in $\A$ is called special precovering if for any object $A\in\A$ there exists an exact sequence $0\to K\to X\to A\to 0$ with $X\in\X$ and $K\in\X^\perp$. Special preenveloping classes can be defined dually. Note that a cotorsion pair $(\X,\Y)$ is complete if and only if $\X$ is special precovering and $\Y$ is special preenveloping.

By  a \textbf{special proper $\X$-resolution}, we mean an $\X$-resolution
$$\cdots \to X_1\to X_0\to A\to 0$$
where each $\Ker(X_i\to X_{i-1})\in \X^\perp$ with  $X_{-1}:=A$.  Note that  any object $A\in\A$ has a special proper $\X$-resolution if and only if $\X$ is a special precovering class.


\begin{prop} \label{bounded resdim} Assume that $(\X,\Y)$ is a complete hereditary  cotorsion pair in $\A$. The following assertions are equivalent for any object $A$ in $\A$ and any integer $n\geq 1$.
	\begin{enumerate}
		\item[(a)]  $\X\textnormal{-}\resdim_\A(A)\leq n.$ 
		\item[(b)] For any special proper $\X$-resolution $$\cdots\to X_n\to X_{n-1}\to\cdots\to X_1\to X_0\to A\to 0$$
		 with $X_{-1}=A$, we have $K_{n-1}:=\Ker(X_{n-1}\to X_{n-2})\in \X$.
		\item[(c)] There exists a special proper $\X$-resolution $$0\to X_n\to\cdots\to X_0\to A\to 0$$
		\item[(d)] $\Ext^{k}_\A(A,Y)=0$ for all objects $Y\in\Y$ and all $k\geq n+1$.
		\item[(e)] $\Ext^{n+1}_\A(A,Y)=0$ for all objects $Y\in\Y$.
	\end{enumerate}
	
	Consequently, $$\X\textnormal{-}\resdim_\A(A)=\sup\{i\in\mathbb{N}: \Ext^{i}_\A(A,Y)\neq 0 \text{ for some } Y\in\Y\}.$$
\end{prop}
\proof $(a)\Rightarrow (b)$ We proceed by induction.

If $n=1$, that is, $\X\textnormal{-}\resdim_\A(A)\leq 1$, then there exists an exact sequence $0\to Z_1\to Z_0\to A\to 0$, with $Z_1,Z_0\in\X$. Let us show that $K_0:=\Ker(X_0\to A)\in \X$. We have the following commutative diagram 
$$\xymatrix{ &0\ar[r]  & Z_1\ar[r]\ar[d]  & Z_0\ar[r]\ar[d] & A\ar[r]\ar@{=}[d] & 0 \\ 
	&0\ar[r]  & K_0\ar[r]  & X_0\ar[r] & A\ar[r]  &0 \\ }$$
which induces the mapping cone that leads to the short exact sequence 
$0\to Z_1\to Z_0\oplus K_0\to X_0\to 0$.
Therefore, $K_0\in \X$ as $Z_1, X_0\in\X$ and $\X$ is closed under extensions and direct summands.

For the case $n\geq 2$, there exists an exact sequence $0\to Z\to Z_0\to A\to 0$, with $Z_0\in\X$ and $\X\textnormal{-}\resdim_\A(Z)\leq n-1$. Since $K_0\in\X^\perp$, we can construct the following commutative diagram with exact rows and columns:
$$\xymatrix{ & & 0\ar[d] & 0\ar[d] & 0\ar[d] &\\ 
	& 0\ar[r] & K_1 \ar[r]\ar[d] & L  \ar[r]\ar[d] & Z \ar[r]\ar[d] & 0\\
	& 0\ar[r] & X_1 \ar[r]\ar[d] & X_1\oplus Z_0  \ar[r]\ar[d] & Z_0 \ar[r]\ar[d] & 0\\
	& 0\ar[r] & K_0 \ar[r]\ar[d] & X_0 \ar[r]\ar[d] & A \ar[r]\ar[d] & 0\\
	& &0  &0 &0 &
}$$

Note that $L\in\X$ as $(\X,\Y)$ is a hereditary cotorsion pair and $X_0,X_1,Z_0\in\X$. Now, $\X\textnormal{-}\resdim_\A(Z)\leq n-1$ and $Z$ has the following special proper $\X$-resolution: $$\cdots \to X'_{n-2}=X_{n-1}\to X'_{n-3}=X_{n-2}\to \cdots \to X'_1=X_2 \to X'_0=L \to Z\to  0$$
By induction, $K_{n-1}=\Ker(X_{n-1}\to X_{n-2})=\Ker(X'_{n-2}\to X'_{n-3})\in\X$ as desired.

$(b)\Rightarrow (c)$ By hypothesis and the fact that $\X$ is a special precovering class.

$(c)\Rightarrow (d)$ Since $(\X,\Y)$ is a hereditary cotorsion pair, for any object $Y\in\Y$ and any integer $k\geq n+1$, we get the following isomorphisms:
$$\Ext^{k}_\A(A,Y)\cong \Ext^{k-1}_\A(K_0,Y)\cong \cdots \cong \Ext^{k-n}_\A(K_{n-1},Y)=\Ext^{k-n }_\A(X_n,Y)=0.$$

$(d)\Rightarrow (e)$ Clear.

$(e)\Rightarrow (a)$ Consider a special proper $\X$-resolution: $$\cdots\to X_n\to\cdots\to X_0\to A\to 0.$$ Since $K_n\in\Y$ and each $X_i\in\X$, $\Ext^{1}_\A(K_{n-1},K_n)\cong \cdots \cong \Ext^{n+1}_\A(A,K_n)=0$. 
Hence the short exact sequence $0\to K_{n}\to X_n\to K_{n-1}\to 0$
splits. Therefore $K_{n-1}\in\X$ and $\X\textnormal{-}\resdim_\A(A)\leq n$.\cqfd

Using the description of the $\X$-resolution dimension given in Proposition \ref{bounded resdim}, the following result is standard and straightforward.
\begin{prop}\label{prop of resdim} Let $(\X,\Y)$ be a complete hereditary cotorsion pair in $\A$.
	\begin{enumerate}
		\item Given a short exact sequence	$\mathcal{E}: 0\to M\to N\to L\to 0$ in $\A$, we have:
		\begin{enumerate}
			\item[(a)] $\X\textnormal{-}\resdim_\A(M)\leq max\{\X\textnormal{-}\resdim_\A(N),\X\textnormal{-}\resdim_\A(L)-1\}.$\\
			The equality holds when $\X\textnormal{-}\resdim_\A(N)\neq \X\textnormal{-}\resdim_\A(L).$
			\item[(b)] $\X\textnormal{-}\resdim_\A(N)\leq max\{\X\textnormal{-}\resdim_\A(M),\X\textnormal{-}\resdim_\A(L)\}.$\\
			The equality holds when $\X\textnormal{-}\resdim_\A(L)\neq \X\textnormal{-}\resdim_\A(M)+1.$
			\item[(c)] $\X\textnormal{-}\resdim_\A(L)\leq max\{\X\textnormal{-}\resdim_\A(N),\X\textnormal{-}\resdim_\A(M)+1\}.$\\
			The equality holds when $\X\textnormal{-}\resdim_\A(M)\neq \X\textnormal{-}\resdim_\A(N).$
		\end{enumerate}
		\item For any family $(M_i)_{i=1,\cdots,n}$ of objects in $\A$, we have  $$\X\textnormal{-}\resdim_\A(\oplus_{i=1}^{n} M_i)=\sup\{\X\textnormal{-}\resdim_\A(M_i)|i=1,\cdots,n\}.$$
	\end{enumerate}
	
	Consequently, $\widehat{\X}$ is thick.
\end{prop}
\begin{lem}\label{core of cot pairs} Let  $\X$ and $\G$ be two classes of objects of $\A$ such that  $\X$ is closed under direct summands. Set $\H=\G^\perp$ and $\Y=\X^\perp$. If $\X$ is a cogenerator for $\G$, then $\G\cap\H\subseteq \X\cap\Y.$
	
	If, in addition, $(\G,\H)$ is a complete cotorsion pair and $\G\cap \widehat{\X}=\X$, then  $$\G\cap\H=\X\cap\Y.$$

\end{lem}
\proof  Let $A$ be an object in $\A$. For the first inclusion, assume that $A\in \G\cap\H$. Since $\X$ is a cogenerator for $\G$, $\X \subseteq  \G$ and then $A\in \H\subseteq \Y$. Moreover, there exists an exact sequence $0\to A\to X\to G\to 0$ with $X\in\X$ and $G\in \G$. This sequence splits as $A\in \H$ and $G\in \G$. Hence, $A\in\X$.

Conversely, assume that that  $A\in \X\cap\Y$. Then, $A\in\X\subseteq \G$. Since the cotorsion pair $(\G,\H)$ is complete, there exists an exact sequence $$0\to A\to W\to G\to 0$$
with $W\in\H$ and $G\in \G$.  Note that $W\in \G\cap\H$ and $\G\cap\H\subseteq \X\cap\Y$ by the previous inclusion. Now, we have $G\in\G$ with $\resdim_\X(G)\leq 1<\infty$, that is, $G\in \G\cap \widehat{\X}$. Hence, $G\in\X$ and this sequence splits. Therefore, $A\in\H$. \cqfd

The following two results relate cotorsion pairs and Hovey triples with (weak) Ab-contexts in a convenient way. The first one is inspired by two results due to  Liang and Yang \cite[Proposition 2.5 and Proposition 2.10]{LY22}.

\begin{prop}\label{AB-cont-cot-pair} Let  $\X$ and $\G$ be two classes of objects of $\A$ such that  $\X$ is closed under direct summands. Set $\Y=\X^\perp$ and $\H=\G^\perp$.
	\begin{enumerate}
		\item 	Assume that $\X$ is a cogenerator for $\G$ and  $\G\cap \widehat{\X}=\X$. If $(\G,\H)$ is a hereditary complete cotorsion pair, then  $(\mathcal{G},\widehat{\X\cap\Y},\X\cap\Y)$ is a left weak AB-context. 
		
		\item Assume that $(\X,\Y)$ is a complete hereditary cotorsion pair. If $(\G,\widehat{\X\cap\Y},\X\cap\Y)$ is a left AB-context, then $(\G,\H)$ is a hereditary complete cotorsion pair with $\H=\widehat{\X\cap\Y}=\widehat{\X}\cap\Y$.
	\end{enumerate}
\end{prop}
\proof 1. By \cite[Proposition 2.5]{LY22}, the pair $(\mathcal{G},\mathcal{G}\cap\mathcal{H})$ is a Frobenius pair in the sense of \cite[Definition 2.5]{BMPS19}. Moreover, using the one-to-one correspondence between left AB-contexts and left Frobenius from \cite[Theorem 5.4(1)]{BMPS19}, we get that the triple $(\mathcal{G},\widehat{(\mathcal{G}\cap\mathcal{H})},\mathcal{G}\cap\widehat{(\mathcal{G}\cap\mathcal{H})})=(\mathcal{G},\widehat{(\mathcal{G}\cap\mathcal{H})},\mathcal{G}\cap\mathcal{H})$ is a left weak AB-context. Finally, applying Lemma \ref{core of cot pairs}, we get our desired left weak AB-context $(\mathcal{G},\widehat{\X\cap\Y},\X\cap\Y)$.

2. By \cite[Theorem 1.12.10]{Has00}, we get that $(\G,\widehat{\X\cap\Y})$ is a complete hereditary cotorsion pair with $\H=\widehat{\X\cap\Y}$. Now let us show that $\H=\widehat{\X}\cap\Y$. 

If $A\in\H$, then there exists a finite $(\X\cap \Y)$-resolution
$$0\to X_n\to\cdots\to X_0\to A\to 0.$$
Since each $X_i\in \X\cap \Y\subseteq \X$, $A\in \widehat{\X}$. We also know that $(\X,\Y)$ is hereditary, then $A\in \Y$. Hence, $A\in \widehat{\X}\cap\Y$. 

Conversely, if $A\in \widehat{\X}\cap\Y$, then $A\in \widehat{\X}$. Moreover, by Proposition \ref{bounded resdim},  there exists a finite special proper $\X$-resolution 
$$0\to X_n\to\cdots\to X_0\to A\to 0$$ 
Since $\Y$ is closed under extensions, each $X_i\in\Y$ and then $X_i\in\X\cap\Y$. Hence, $A\in \widehat{\X\cap\Y}=\H$. \cqfd

\bigskip
Assume that $(\Q,\W,\R)$ is  a hereditary Hovey triple  and let
$$(\widetilde{\Q},\R):=(\Q\cap \W,\R) \text{ and } (\Q,\widetilde{\R}):=(\Q, \W\cap\R)$$
be the associated cotorsion pairs. Then, the  following two assertions always hold:

$\bullet$  $(\widetilde{\Q},\R)$ is a complete hereditary cotorsion pair.

$\bullet$ $\left(\Q,\widehat{(\Q\cap\widetilde{\R})},\widetilde{\Q}\cap\R\right)$ is a left weak AB-context: indeed, since $(\Q,\widetilde{\R})$ is a complete hereditary cotorsion pair and $\Q\cap\widetilde{\R}=\widetilde{\Q}\cap\R$,  we get this weak AB-context by letting $\G=\X=\Q$ in Proposition \ref{AB-cont-cot-pair}(1).

Under the finiteness assumption $\A=\widehat{\G}$, the following result shows that there is a converse to this, giving a new characterization of hereditary Hovey triples in terms of AB-contexts. We note that such a relation between Hovey triples and (weak) AB-contexts was also noted by A. Xu in \cite[Theorem  4.2]{Xu17} in a particular setting.

\begin{thm}\label{Hov from AB-cont} Let $\X\subseteq \G\subseteq\A$ be two classes such that:
	\begin{enumerate}
		\item $(\X,\Y)$ is a complete hereditary  cotorsion pair.
		\item $(\G,\H,\X\cap\Y)$ is a left AB-context.
	\end{enumerate}
	
	Then,  $(\G,\widehat{\X},\Y)$ is a hereditary Hovey triple.
\end{thm}
\proof We divide the proof into three parts:

(a)  The class  $\widehat{\X}$ is thick  by Proposition \ref{prop of resdim}.

(b) By assumption and \cite[Theorem 1.12.10(1)]{Has00}, $\H=\widehat{\X\cap \Y}$. Then, $\H=\widehat{\X}\cap \Y$ and  the pair $(\G,\widehat{\X}\cap \Y)$ is a hereditary complete cotorsion pair by  Proposition \ref{AB-cont-cot-pair}(2).

(c) It remains to show that $(\G\cap\widehat{\X},\Y)$ is a hereditary complete cotorsion pair. By assumptions, we only need to show  $\X=\G\cap \widehat{\X}$. 

The inclusion $\X\subseteq \G\cap \widehat{\X}$ is clear. 

Conversely, assume $A\in\G\cap \widehat{\X}$. Then, $n=\X\textnormal{-}\resdim(A)<\infty$. If $n=0$, that is, $X\in\X$, then we are done. Assume now that $n\geq 1$. Let us proceed by induction on $n$. 

There exists an exact sequence 
$0\to K\to X\to A\to 0$
with
$\X\textnormal{-}\resdim_\A(K)=n-1$ and $X\in \X$. Since $(\X,\Y)$ is complete, there exists an exact sequence 
$0\to K\to Y\to L\to 0$
with $\Y\in\Y$ and $L\in\X$.  Consider now the following pushout:
$$\xymatrix{ & 0\ar[d] & 0\ar[d] & & \\  0\ar[r] & K\ar[d] \ar[r] & X\ar@{-->}[d] \ar[r] & A \ar@{=}[d] \ar[r] & 0 \\ 0\ar[r] & Y \ar@{-->}[r]\ar[d] & H \ar[d] \ar[r] & A\ar[r] & 0 \\ & L\ar@{=}[r] \ar[d] & L\ar[d] & & \\ & 0 & 0}.$$
Since $A\in\G$, $X\in\X\subseteq \G$ and $\G$ is closed under kernels of epimorphisms, $K\in \G$. So, by induction, $K\in\X$. Hence, $Y\in\X\cap\Y$. 

On the other hand, $\X\cap\Y=\G\cap\H$ by hypothesis. Then, the middle short exact sequence splits as $A\in \G$ and $Y\in \H=\G^\perp$. But, $H\in \X$ by the middle column. Thus, $A\in\X$ as desired. \cqfd

We are now ready to prove our first main result of this section (Theorem D from the Introduction). We only state and prove the first assertion as the second one is completely dual.

\begin{thm}\label{const Hov trip} Assume that $(\Q,\widetilde{\R})$ and $(\widetilde{\Q},\R)$ are two complete hereditary cotorsion pairs such that:
	\begin{enumerate}
		\item[(b)] $\widetilde{\Q}\subseteq \Q$ (or equivalently, $\widetilde{\R}\subseteq\R $).
		\item[(c)] $\Q\cap \widetilde{\R}=\widetilde{\Q}\cap \R$.
	\end{enumerate}
	
	If $\Q\textnormal{-}\resdim_{\A}(\A)<\infty$, then $(\Q,\widehat{(\widetilde{\Q})},\R)$ is a hereditary Hovey triple.
	
\end{thm}
\proof Taking advantage of our previous notations, let us set  $$(\G,\H)=(\Q,\widetilde{\R}) \text{ and }(\X,\Y)=(\widetilde{\Q},\R).$$ By Theorem \ref{Hov from AB-cont} and hypothesis, it suffices to show that
$$(\G,\H,\X\cap\Y)=(\Q,\widetilde{\R},\widetilde{\Q}\cap\R)$$ is a left weak AB-context.  

Notice that $\G$ is a cogenerator for $\G$ with $\G\cap \widehat{\G}=\G$, so $(\G,\widehat{\G\cap\H},\G\cap\H)$  is a left AB-context by Proposition \ref{AB-cont-cot-pair}(1) and so a left AB-context as $\A=\widehat{\Q}$. Then, since $\G\cap\H=\X\cap\Y$, we can apply Proposition \ref{AB-cont-cot-pair}(2) to get that $\H=\widehat{\X\cap\Y}$ and so that $(\G,\H,\X\cap\Y)$ is a left AB-context as desired. \cqfd

\bigskip

Now we turn our attention to our last main result (Theorem C from the Introduction). 

In \cite{BEGO22a, BEGO22b} the notions of $\F$-flat and $\C$-cotorsion modules were introduced as relative versions of flat and cotorsion modules, respectively. These two notions preserve most of the properties satisfied by the flat and cotorsion modules. But some properties related to our $\GC$-flat model structure escape this extension. For instance, the pair $(\F(R),\C(R))$ is not a cotorsion pair unless $R\in\F(R)$ \cite[Theorem 3.7(3)]{BEGO22b}. However, it turns out that there are two other classes of modules,
$$\V_C(R)=\;^\perp(\B_C(R)\cap\C(R)) \text{ and }\H_C(R)=\B_C(R)\cap\C(R),$$
which could serve as relative versions of the class of flat and cotorsion modules, respectively. Under some assumptions (depending on each property), here are some of their properties that support this claim:

\begin{enumerate}
	\item[(a)] If $C=R$, then $\F(R)=\V_C(R)=\Flat(R) \text{ and }\C(R)=\H_C(R)=\Cot(R).$
	\item[(b)] $(\V_C(R),\H_C(R))$ is a complete hereditary cotorsion pair.
	\item[(c)] $\F(R)\subseteq \V_C(R)\subseteq \GCF(R).$
	\item[(d)] $\GCC(R)\subseteq \H_C(R)\subseteq \C(R).$
	\item[(e)] $\F(R)\cap\C(R)=\V_C(R)\cap\H_C(R)=\GCF(R)\cap\GCC(R).$
\end{enumerate}

Based on this discussion, it is natural to introduce the following definition. 
\begin{defn}An $R$-module $H$ is said to be $\H_C$-\textbf{cotorsion} if it belongs to $\H_C(R):=\B_C(R)\cap \C(R)$. An $R$-module $V$ is said to be $\V_C$-\textbf{flat} if $\Ext^1_R(V,H)=0$ for all $\H_C$-cotorsion  $R$-modules.
	
	Set $\V_C(R):=\;^\perp\H_C(R)$ the class of all $\V_C$-flat $R$-modules.
	
	Given an $R$-module $M$, we define its \textbf{$\V_C$-flat  dimension} as $$\V_C\textnormal{-}\fd_R(R)=\V_C(R)\textnormal{-}\resdim_R(M).$$
\end{defn}

We recall one last ingredient and we are ready to prove Theorem C. An $R$-module $C$ is w-tilting \cite[Definition 2.1]{BGO16a} if $\Ext_R^{i\geq 1}(C,C^{(I)})=0$ for every set $I$ and $_RR$ has a $\Hom_R(-,\Add_R(C))$-exact $\Add_R(C)$-coresolution. From \cite[Section 2]{BGO16a} and \cite[Proposition 2.14]{BEGO22b}, we have the relationships
$(a)\Rightarrow (b) \Rightarrow (c)$,
among the following assertions:

(a) $C$ is a semidualizing $(R,S)$-bimodule. 

(b) $_RC$ and $C_S$ are $w$-tilting having a degreewise finite projective resolution.

(c) $_RC$ and $C_S$ are $w^+$-tilting modules.


\begin{thm}\label{GCF model structure} Assume that $R$ is $\GCF$-closed and $C$ is $w$-tilting admitting a degreewise finite projective resolution.  If $R$ has finite global $\GC$-flat dimension, then there exists an abelian model structure on $R$-Mod,  $$\left( \GCF(R),\widehat{\V_C(R)},\H_C(R)\right),$$ as follows:
	
	$\bullet$ The cofibrant objects coincide with $\GC$-flat $R$-modules.
	
	$\bullet$ The trivial objects coincide with modules having finite $\V_C$-flat dimension.
	
	$\bullet$ The fibrant objects coincide with $\H_C$-cotorsion $R$-modules.
	
\end{thm}
\proof  We know by \cite[Theorem 4.15, Proposition 4.13 and Proposition 4.14]{BEGO22b}, that
$$(\GCF(R),\GC(R)) \text{ and } (\V_C(R),\H_C(R))$$ are two complete hereditary cotorsion pairs with with the same core, that is, $$\GCF(R)\cap\GCC(R)=\V_C(R)\cap\H_C(R).$$
Moreover, one can see from the proof of \cite[Theorem 4.15]{BEGO22b} that $\GCC(R)\subseteq \H_C(R)$. 
Thus, this result follows by Theorem \ref{const Hov trip}. \cqfd

The following consequence has been proven by A. Xu \cite[Corollary 4.6(3)]{Xu17} when the ring is right coherent.
\begin{cor}For any ring $R$ with finite weak global Gorenstein dimension, we have a hereditary abelian model structure on $R$-Mod  $$(\GF(R),\widehat{\Flat(R)},\Cot(R)).$$
\end{cor}

\bigskip

	\noindent\textbf{Acknowledgement.}

The authors would like to thank the referee for the insightful suggestions that enhanced the exposition.

\end{document}